\magnification=\magstephalf
\input amstex
\loadbold
\documentstyle{amsppt}
\refstyle{A}
\NoBlackBoxes

\vsize=7.5in

\def\pf{\hfill $\square$}
\def\c{\cite}
\def\fb{\frak{b}}

\def\fg{\frak{g}}

\def\fk{\frak{k}}

\def\end{\text{End}}

\def\ds{\Bbb D}
\def\bds{\partial{\Bbb D}}

\topmatter
\title Some remarks on CMV matrices and dressing orbits\endtitle
\leftheadtext{L.-C. Li}
\rightheadtext{Some remarks on CMV matrices and dressing orbits}

\author Luen-Chau Li\endauthor
\address{L.-C. Li, Department of Mathematics,Pennsylvania State University
University Park, PA  16802, USA}\endaddress
\email luenli\@math.psu.edu\endemail
\abstract  The CMV matrices are the unitary analogs of  Jacobi matrices.
In the finite case, it is well-known that the set of Jacobi matrices with a 
fixed trace is nothing but a coadjoint orbit of the lower triangular
group.  In this note, we will give the analog of this result for the
CMV matrices. En route, we also discuss the Hamiltonian formulation
of the Lax equations for the defocusing Ablowitz-Ladik hierarchy.
\endabstract
\endtopmatter

\document
\subhead
1. \ Introduction
\endsubhead

\baselineskip 15pt
\bigskip

A major development in the theory of orthogonal polynomials on
the unit circle (OPUC) is the introduction of the so-called
CMV matrices by Cantero, Moral and Val\'azquez in 2003 \c{CMV}.
CMV matrices are the unitary analogs of Jacobi matrices (which arise in the 
theory of orthogonal polynomials on the real line (OPRL))and
as clearly demonstrated in the two-volume monograph of Simon \c{S},
they provide a powerful tool in OPUC.  It is well-known that
Jacobi matrices are associated with one of the most
celebrated examples of integrable systems, namely, the Toda
lattice and much has been written about the subject.(See,
for example, the monograph \c{T} and the references therein.)
In view of this, it is natural to ask if there is an integrable
system related to the CMV matrices and OPUC in the same way
that Toda relates to Jacobi matrices and OPRL.  The answer
to this question was obtained in the recent work
of Nenciu \c{N1},\c{N2}, who showed that the sought-after integrable 
system is the defocusing Ablowitz-Ladik (AL) system (a.k.a. defocusing 
discrete nonlinear Schr\"odinger equation) \c{AL}.  More specifically,
by making use of the connection between Ablowitz-Ladik and
OPUC, the author in \c{N1},\c{N2} has successfully derived
the Lax pair formulation for the nonlinear equation.

Returning to the Jacobi operators, it is well-known that in the finite
case, the collection of such operators with a fixed trace is
a coadjoint orbit of the lower triangular group, when the dual
of its Lie algebra is identified with the real symmetric matrices. 
(See, for example, \c{A}, \c{K} and \c{R}.)  In the case of the finite 
CMV matrices, it is  natural to ask if there is an analogous Poisson 
geometric meaning.  In this short note, we will provide an answer to this
question and en route, we will discuss the Hamiltonian formulation
of the Lax equations for the defocusing AL hierarchy.  As the reader
will see in Section 2, we have a coboundary Poisson Lie group
$(G^{\Bbb R}, \{\cdot,\cdot\}_{J})$ (in the sense of Drinfe'ld \c{D})
containing the unitary group $U(n)$ as a Poisson Lie subgroup.  Also, 
there exists a special CMV matrix $x_{f}$ with 
$\theta$-factorization $x_{f}= x^{e}_{f} x^{o}_{f}$ (we are
following the terminology of \c{S} here). If 
${\Cal L}_{x^{e}_{f}}$ and ${\Cal L}_{x^{o}_{f}}$ are
the dressing orbits of the dual group $G^{\Bbb R}_{J}$ through  $x^{e}_{f}$ 
and 
$x^{o}_{f}$ respectively, then our main result is the following:
the factors $g^{e}$ and $g^{0}$ in the (unique)
$\theta$-factorization of a CMV matrix $g({\underline\alpha})$ are 
respectively  elements of ${\Cal L}_{x^{e}_{f}}$ and ${\Cal L}_{x^{o}_{f}}$;
moreover, the set of CMV matrices is the image of the symplectic leaf
${\Cal L}_{x^{e}_{f}}\times {\Cal L} _{x^{o}_{f}}$ of $U(n)\times U(n)$
(equipped with the product structure) under the Poisson
automorphism
$ m\mid {\Cal L}_{x^{e}_{f}}\times {\Cal L} _{x^{o}_{f}}:
     {\Cal L}_{x^{e}_{f}}\times {\Cal L} _{x^{o}_{f}}
     \longrightarrow \{\hbox{CMV}\,\,\hbox{matrices}\}$
where $m:U(n)\times U(n)\longrightarrow U(n)$ is the multiplication
map of the group $U(n).$

\bigskip
\noindent{\bf Acknowledgments.} The author would like to thank Percy
Deift for showing him the thesis of Nenciu  and for pointing out that 
the CMV matrices are unitary.  He is also grateful to Irina Nenciu for 
sending him a copy of her thesis.
\bigskip
\bigskip

\subhead
2. \ CMV matrices and dressing orbits
\endsubhead

\bigskip

Let $G^{\Bbb R}$ be $GL(n,\Bbb C)$ considered as a real
Lie group, and let $K$ and $B$ be respectively the
unitary group $U(n)$ and the lower triangular group
with positive diagonal entries.  It is well-known that
$G^{\Bbb R}$ admits the Iwasawa decomposition
$$G^{\Bbb R} = KB \eqno(2.1)$$
which means every $g\in G^{\Bbb R}$ admits a factorization
$$g = g_{+} g _{-}^{-1}\eqno(2.2)$$
for unique $g_{+}\in K$ and $g_{-}\in B.$(We shall henceforth
use the notation $g_{+}$ and $g_{-}$ with this interpretation.) On the
Lie algebra level, we have 
$$\fg^{\Bbb R} = \frak k \oplus \frak b, \eqno(2.3)$$
where $\fg^{\Bbb R}$, $\fk$ and $\fb$ are respectively
the Lie algebras of $G^{\Bbb R}$, $K$ and $B$.
We shall denote by $\Pi_{\frak k}$ and $\Pi_{\frak b}$ the
projection maps onto $\frak k$ and $\frak b$ relative to the splitting 
in (2.3) and set
$$J = \Pi_{\frak k} -\Pi_{\frak b}.\eqno(2.4)$$
From classical r-matrix theory \c{STS1},\c{STS2},
$J:\fg^{\Bbb R}\longrightarrow \fg^{\Bbb R}$ is a solution of the modified 
Yang-Baxter equation (mYBE).  Consequently, we can equip
$\fg^{\Bbb R}$ with the $J$-bracket
$$[X,Y]_{J} = {1\over 2}([JX,Y] + [X,JY]).\eqno(2.5)$$
In what follows, we shall denote the vector space $\fg^{\Bbb R}$
equipped with the $J$-bracket by $\fg^{\Bbb R}_J$. Indeed, it
is easy to check from (2.5) that $\fg^{\Bbb R}_J= \fk\ominus\fb$
(Lie algebra anti-direct sum).

Note that explicitly, the projection maps 
are given by the formulas
$$\Pi_{\frak k} X = X_{+}- (X_{+})^*,\quad 
  \Pi_{\frak b} X = X_{-} + X_{0} + (X_{+})^*,\eqno(2.6)$$
where $X_{+}$, $X_{0}$ and $X_{-}$
are the strict upper, diagonal and strict lower
parts of $X\in \fg^{\Bbb R}$ and we will make use of these in
the sequel.  As non-degenerate invariant pairing on $\fg^{\Bbb R}$, we 
take $$(X,Y) = Im\, tr (XY).\eqno(2.7)$$
This choice is critical for what we have in mind and with respect
to $(\cdot,\cdot)$, we  now define the right and left gradients
of a smooth function $\varphi$ on $G^{\Bbb R}$ by
$$(D\varphi(g), X)= {d\over dt}{\Big|_{t=0}} \varphi(e^{tX}g),
\quad (D' \varphi(g), X)= {d\over dt}{\Big|_{t=0}} 
\varphi(ge^{tX}), \, X\in \fg^{\Bbb R}.\eqno(2.8)$$

\proclaim
{Proposition 2.1} (a) $J$ is skew-symmetric relative to the
pairing $(\cdot,\cdot).$  Consequently, $(G^{\Bbb R}, \{\cdot,\cdot\}_{J})$ is 
a coboundary Poisson Lie group with tangent Lie bialgebra 
$(\fg^{\Bbb R},\fg^{\Bbb R}_J)$
where $\{\cdot,\cdot\}_{J}$ is the Sklyanin bracket
$$\{\varphi, \psi\}_{J}(g) =(J(D'\varphi(g)), D'\psi(g)) - 
(J(D\varphi(g)), D\psi(g)).\eqno(2.9)$$
Moreover, $K=U(n)$ is a Poisson Lie subgroup of 
$(G^{\Bbb R}, \{\cdot,\cdot\}_{J}).$ 
\newline
(b) The Hamiltonian equation of motion generated by 
$\varphi\in C^{\infty}(G^{\Bbb R})$ 
is given by
$$\dot g = g (\Pi_{\fk}(D'\varphi(g)))-(\Pi_{\fk}(D\varphi(g)))g.\eqno(2.10)$$
In particular, for the Hamiltonian $H_{k}(g) ={1\over k} Re\, tr\, g^k$, the 
corresponding equation of motion is
$$\dot g = g\left(ig^{k}_{+} +i(g^{k}_{+})^*\right)                                       -\left(ig^{k}_{+} +i(g^{k}_{+})^*\right)g.\eqno(2.11)$$
\newline
(c) The underlying group of the Poisson group $G^{\Bbb R}_{J}$ dual to
$(G^{\Bbb R}, \{\cdot,\cdot\}_{J})$ consists of $G$ equipped with
the multiplication
$$g\ast h \equiv g_{+}hg_{-}^{-1}.\eqno(2.12)$$
\newline
(d) Equip $G^{\Bbb R}\times G^{\Bbb R}$ with the product structure.
Then the Hamiltonian equations of motion generated by
${\widetilde H}_{k}(g_1,g_2) ={1\over k} Re\, tr\,(g_{1}g_{2})^{k}$
are given by
$$\eqalign{
\dot g_1 &= g_{1} (\Pi_{\fk}(i(g_2 g_1)^k)) -(\Pi_{\fk} (i(g_1 g_2)^k))g_1, \cr
\dot g_2 &= g_{2}(\Pi_{\fk} (i(g_1 g_2)^k)) -(\Pi_{\fk} (i(g_2 g_1)^k))g_2.\cr}
\eqno(2.13)
$$
Moreover, the monodromy matrix $g= g_1g_2$ satisfies (2.11).
\endproclaim 

\demo
{Proof} (a) Since $\fk$ is a real form of $gl(n,\Bbb C)$, it
follows that $tr (XY)\in \Bbb R$ for $X,Y\in \fk.$  Consequently,
$\fk$ is an isotropic subalgebra of $\fg^{\Bbb R}$ relative to
$(\cdot,\cdot)$, i.e., $(\fk,\fk)=0$.
On the other hand, $\fb$ is also an isotropic subalgebra of 
$\fg^{\Bbb R}$ relative to $(\cdot,\cdot)$ because the diagonal
entries of the elements in $\fb$ are real.  Combining these
two facts, it follows that $J$ is skew-symmetric relative
to $(\cdot,\cdot)$.  The rest of the assertion concerning
$(G^{\Bbb R}, \{\cdot,\cdot\}_{J})$  then follows
from standard results in \c{STS2}.  Finally, in order to
show that $K$ is a Poisson Lie subgroup, it suffices to
check that $K$ is a Poisson submanifold of 
$(G^{\Bbb R}, \{\cdot,\cdot\}_{J}).$  We shall leave
the simple verification to the reader.
\newline
\noindent (b) The calculation is standard.  Note that in 
deriving (2.11), we have made use of the formula
$D'H_{k}(g)=DH_{k}(g)=ig^{k}$ and the explicit expression
for $\Pi_{\fk}$ in (2.6).
\newline\
\noindent (c), (d) We shall leave the verification to 
the reader.
\pf
\enddemo

The equations in (2.11) together with the fact that $K=U(n)$ is a 
Poisson Lie subgroup of  $(G^{\Bbb R}, \{\cdot,\cdot\}_{J})$ means
that the restriction of these equations to $K$ are Hamiltonian
with respect to the induced structure on $K$.  Moreover, eqn(4.13)
in \c{N1} is a special case of (2.11) above if we take $g$ to
be a finite CMV matrix.  Since the CMV matrix is a very special
unitary matrix (see Definition 2.2 below), we ask the question if 
the collection of such matrices has a
natural Poisson geometric meaning which would allow such a
restriction to happen. Before we turn to answer this question, let us 
recall the definition of a finite CMV matrix \c{N1},\c{N2} in a form which 
is suitable for our purpose here.
To simplify the language, we shall drop the term ``finite'' from
now on.

We begin with some notations.  Let $\ds$ be the open unit disk
$\{z\in {\Bbb C}\mid |z|<1\}$ and let $\bds$ be its boundary.
Given an (n-1)-tuple 
${\underline\alpha} = (\alpha_0,\cdots,\alpha_{n-2})\in \ds^{n-1},$
we define unitary matrices 
$$\eqalign{
\theta_j = \pmatrix
  \bar \alpha_j& \rho_j\cr
  \rho_j& -\alpha_j\cr
  \endpmatrix , \quad \rho_j= (1-|\alpha_j|^2)^{1\over 2},\,\,
  j=0,\cdots,n-2}\eqno(2.14)$$
$$\theta_{n-1} = -1.\eqno(2.15)$$

\definition
{Definition 2.2} The CMV matrix associated with an (n-1)-tuple
${\underline\alpha} = (\alpha_0,\cdots,\alpha_{n-2})\in \ds^{n-1}$
is the penta-diagonal unitary matrix given by
$$g({\underline\alpha}) = g^{e}\left(\{\alpha_{2j}\}_{j=0}^{[\frac{n-2}{2}]}\right)
g^{0}\left(\{\alpha_{2j+1}\}_{j=0}^{[\frac{n-3}{2}]}\right)\eqno(2.16)$$
where
$$g^{e}\left(\{\alpha_{2j}\}_{j=0}^{[\frac{n-2}{2}]}\right) =
diag\left(\theta_0,\theta_2,\cdots,\theta_{2[\frac{n-1}{2}]}
\right)
\eqno(2.17)$$
and
$$g^{0}\left(\{\alpha_{2j+1}\}_{j=0}^{[\frac{n-3}{2}]}\right) = 
diag\left(1,\theta_1,\theta_3,\cdots,\theta_{2[\frac{n-2}{2}]
+1}\right).\eqno(2.18)$$
\enddefinition
\medskip
\noindent{\bf Remark 2.3} (a) It follows from the above definition that the
$(2j, 2j+2)$ and the $(2j+1, 2j-1)$ entries of a CMV matrix are
zero for any $j$. The factorization in (2.16) above is called
the $\theta$-factorization following the terminology in \c{S}.
We have more to say on this below.
\smallskip
\noindent (b) In the original definition of the CMV matrix in \c{N1},\c{N2},
there is an extra parameter $\alpha_{n-1}\in \bds$ involved.  But
subsequently, the author restricts her attention to $\alpha_{n-1} = -1.$ 
The reader will see that this assumption is natural from our point of
view.
\smallskip
\noindent (c) Given a nontrivial probability measure $d\mu$ on
$\bds$, Cantero, Moral and Val\'azquez \c{CMV} produces an
orthonormal basis of $L^{2}(\bds, d\mu)$ by applying Gram-Schmidt
to $1,z,z^{-1},z^{2},z^{-2},\cdots.$  The matrix representation of
the operator $f(z)\mapsto zf(z)$ in  $L^{2}(\bds, d\mu)$ in this
basis is an infinite CMV matrix.  The finite case which we consider here
corresponds to a trivial probability measure $d\mu$ supported
at $n$ points and the $\alpha_j$'s are the Verblunsky coefficients
which appear in Szeg\H o recursion \c{S}.
\medskip

In general, we shall denote by $g^{e}$ any $n\times n$ block
diagonal matrix with $2\times 2$ blocks on the main diagonal of
the form
$$\eqalign{
\ \pmatrix
  \bar \alpha& \rho\cr
  \rho& -\alpha\cr
  \endpmatrix , \quad \rho = (1-|\alpha|^2)^{1\over 2},\quad \alpha\in \ds \cr}
\eqno(2.19)$$                                           
except when $n$ is odd, the last block is the number $-1.$
We shall denote the collection of such matrices by ${\Cal T}^{e}$.
Similary, we shall denote by $g^{o}$ any $n\times n$ block
diagonal matrix which begins with the $1\times 1$ block equal
to $1$ followed by $2\times 2$ blocks of the form in (2.19)
except when $n$ is even, the last block is the number $-1.$
We shall use the symbol ${\Cal T}^{o}$ to denote the
collection of such matrices.  Clearly, for given $g^{e}\in {\Cal T}^{e}$
and $g^{o}\in {\Cal T}^{o}$, there exists unique 
${\underline\alpha} = (\alpha_0,\cdots,\alpha_{n-2})\in \ds^{n-1}$
such that
$$g^{e}g^{o} = g({\underline\alpha}).\eqno(2.20)$$
Indeed, more is true, namely, it is straightforward to verify
that the map
$$\eqalign{
   m\mid {\Cal T}^{e}\times {\Cal T}^{o}:&{\Cal T}^{e}\times {\Cal T}^{o}
   \longrightarrow \{\hbox{CMV}\,\,\hbox{matrices}\}\cr
    &(g^{e},g^{o})\mapsto g^{e}g^{o}\cr}\eqno(2.21)$$
is a diffeomorphism, where $m:K\times K\longrightarrow K$
is the multiplication map of the group $K$.

In order to understand the Poisson geometric meaning of the 
collection of CMV matrices,
we  appeal to the following result in the theory of Poisson Lie groups:
the symplectic leaves of a Poisson Lie group are given by the orbits of 
so-called dressing actions \c{STS2},\c{LW}.  Indeed, it follows from 
Theorem 13 of \c{STS2} (which applies to the coboundary case)
that the  symplectic leaf of  $(G^{\Bbb R}, \{\cdot,\cdot\}_{J})$   
passing through $x\in G^{\Bbb R}$ is given
by
$${\Cal L}_{x} =\left\{\, g_{+}^{-1}x(x^{-1}gx)_{+}\mid g\in G^{\Bbb R}_J \,
\right\}.
\eqno(2.22)$$

In analogy with Example 4.2.7 in \c{S}, we introduce the following
special CMV matrix
$$x_{f} = x^{e}_{f} x^{o}_{f} \eqno(2.23)$$
corresponding to 
$${\underline\alpha} = (0,0,\cdots, 0).\eqno(2.24)$$
In other words,
$$x^{e}_{f}= diag(w^*,w^*,\cdots)\eqno(2.25)$$
and
$$x^{o}_{f}= diag(1, w^*,w^*,\cdots)\eqno(2.26)$$
where
$$\eqalign{w^*= \pmatrix
  0&1\cr
  1&0
\endpmatrix .}\eqno(2.27)$$

We now come to the main result of this work.

\proclaim
{Theorem 2.4} (a) ${\Cal L}_{x^{e}_{f}} = {\Cal T}^{e}.$
\newline
(b) ${\Cal L} _{x^{o}_{f}} = {\Cal T}^{o}.$
\newline
(c) The product ${\Cal L}_{x^{e}_{f}}\times {\Cal L} _{x^{o}_{f}}$ is 
a symplectic leaf of $K\times K$ equipped with the product structure.
Moreover, the collection of CMV matrices is the image of 
${\Cal L}_{x^{e}_{f}}\times {\Cal L} _{x^{o}_{f}}$ under the 
Poisson automorphism 
$ m\mid {\Cal L}_{x^{e}_{f}}\times {\Cal L} _{x^{o}_{f}}:
     {\Cal L}_{x^{e}_{f}}\times {\Cal L} _{x^{o}_{f}}
     \longrightarrow \{\hbox{CMV}\,\,\hbox{matrices}\}$
where $m:K\times K\longrightarrow K$ is the multiplication
map of the group $K.$
\endproclaim

\demo
{Proof} (a) Take an arbitrary element
$$\eqalign{a &= g_{+}^{-1}x^{e}_{f}((x^{e}_{f})^{-1}gx^{e}_{f})_{+}\cr 
             &=  g_{-}^{-1}x^{e}_{f}((x^{e}_{f})^{-1}gx^{e}_{f})_{-}\cr}\eqno(2.28)$$
in the dressing orbit through $x^{e}_{f}$.  We first consider the
case where $n$ is even. From the first line of
(2.28), it is clear that $a$ is unitary.  On the other hand,
since $g_{-}$ is lower
triangular, it follows from the second line of (2.28) that $a$ is block 
lower triangular with $2 \times 2$ blocks on the diagonal.
Moreover, from the fact that the diagonal entries of $g_{-}$ are
positive, it is easy to see that each of the $2\times 2$ blocks 
on the main diagonal has the following
properties: (i) the entry in the upper right hand corner is positive, 
(ii) the determinant is negative (since $det(w^*) = -1$).  Consequently,
the matrix $(a^{*})^{-1}$ is upper block triangular with diagonal blocks
having the same properties.  But $a = (a^{*})^{-1}$, so it follows
that $a$ must be block diagonal, i.e.,
$$a = diag\left(\phi_0,\phi_2,\cdots,\phi_{2[\frac{n-1}{2}]}
\right)
\eqno(2.29)$$
where for each $j$, $\phi_{2j}$ is a unitary $2\times 2$ matrix 
with a positive entry in the upper right hand corner and whose
determinant is $-1$.  Consequently, $\phi_{2j}$ must be
of the form
$$\eqalign{
\phi_{2j} = \pmatrix
  \bar \alpha_{2j}& \rho_{2j}\cr
  \rho_{2j}& -\alpha_{2j}\cr
  \endpmatrix}\eqno(2.30)$$
for some $\alpha_{2j}\in \ds$, where $\rho_{2j} =(1-|\alpha_{2j}|^{2})^{1\over2}$. 
Hence we have shown that
${\Cal L}_{x^{e}_{f}} \subset {\Cal T}^{e}.$
Conversely, take an arbitrary element
$$g^{e} = diag\left(\theta_0,\theta_2,\cdots,\theta_{2[\frac{n-1}{2}]}
\right)
\eqno(2.31)$$ 
in ${\Cal T}^{e}$ where $\theta_{2j}$ is of the form given in
(2.14).  Define
a block diagonal matrix
$$g= diag\left(l_0,l_2,\cdots,l_{2[\frac{n-1}{2}]}\right)\eqno(2.32)$$
such that
$$\eqalign{
l_{2j} =  \pmatrix
          \rho_{2j}& 0\cr
          -\alpha_{2j}& 1\cr
          \endpmatrix ,\quad j=0,\cdots, \left[\frac{n-1}{2}\right].\cr}
\eqno(2.33)$$
Clearly, $g$ is lower triangular so that
$$\eqalign{&
g_{+}^{-1}x^{e}_{f}((x^{e}_{f})^{-1}gx^{e}_{f})_{+}\cr 
=\, & (gx^{e}_{f})_{+}\cr
=\, & \left(diag(l_{0}w^*,l_{2}w^*,\cdots)\right)_{+}.\cr}
\eqno(2.34)$$
But from the definition of $l_{2j}$, we find that
$l_{2j}w^*$ admits the factorization
$$l_{2j}w^* = \theta_{2j}\pmatrix
\rho_{2j} & 0\cr
-\bar\alpha_{2j}&1\endpmatrix .\eqno(2.35)$$
Hence it follows that 
$g_{+}^{-1}x^{e}_{f}((x^{e}_{f})^{-1}gx^{e}_{f})_{+}= g^{e}.$
Consequently, we have the reverse inclusion
${\Cal T}^{e}\subset {\Cal L}_{x^{e}_{f}}$
as well. When $n$ is odd, everything goes through the same as before 
except that for each of the matrices in (2.29),(2.31)-(2.32),
the last block  
is a $1\times 1$ block.  We shall leave the easy detail to
the reader. 
\newline
\noindent (b) The argument is similar to (a).
\newline
\noindent (c) This is clear from the definition of the CMV matrices.
\pf 
\enddemo

\proclaim
{Corollary 2.5} Let $g^{e}$, $g^{o}$ have their
usual meaning and let $g({\underline\alpha}) = g^{e}g^{0}$
and $\widetilde g({\underline\alpha})= g^{0}g^{e}$ where $\alpha\in \ds$
is uniquely determined by $g^{e}$, $g^{o}$.  Then
the equations
$$\eqalign{
\dot g^{e} &= g^{e} (\Pi_{\fk}(i{\widetilde g({\underline\alpha})}^k)) -
(\Pi_{\fk} (ig({\underline\alpha})^k))g^{e}, \cr
\dot g^{o} &= g^{o} (\Pi_{\fk} (ig({\underline\alpha})^k)) -
(\Pi_{\fk} (i{\widetilde g({\underline\alpha})}^k))g^{o}\cr}
\eqno(2.36)
$$
are the Hamiltonian equations of motion on the symplectic manifold
${\Cal L}_{x^{e}_{f}}\times {\Cal L} _{x^{o}_{f}}$ generated
by the Hamiltonian
${\widetilde H}_{k}(g^{e},g^{o}) ={1\over k} 
Re\, tr\,(g({\underline\alpha}))^{k}.$  Moreover, under the Hamiltonian flow
defined by (2.36), $g({\underline\alpha})$ 
evolves according to
$$\dot g({\underline\alpha}) = g({\underline\alpha})(\Pi_{\fk}
(ig({\underline\alpha})^{k}))-
(\Pi_{\fk}(ig({\underline\alpha})^{k}))g({\underline\alpha}).\eqno(2.37)$$
\endproclaim

\demo
{Proof} This is a consequence of the above theorem and Proposition
2.1 (d).
\pf
\enddemo

\noindent {\bf Remark 2.6} (a) It follows from the r-matrix formulation
that the equations (2.36)-(2.37) can be solved via factorization
problems.  Actually, the same remark also holds true for the infinite
case in \c{N1},\c{N2} as one can extend the Iwasawa decomposition
to the group of bounded invertible operators on ${l}^{2}({\Bbb Z}_+)$
(cf. \c{DLT}).
\newline
\noindent (b) The above corollary suggests that in a sense, it seems more
natural to consider (2.36).  Whether this is so from the point of
view of OPUC remains to be seen.


\bigskip


\newpage

\Refs
\widestnumber\key{STS1}

\ref\key{A}
\by Adler, M.
\paper On a trace functional functional for formal 
pseudodifferential operators and the symplectic structure for 
Korteweg-de Vries type equations
\jour Invent. Math. \vol 50\yr 1979\pages 219-248
\endref

\ref\key{AL}
\by Ablowitz, M. and Ladik, J.
\paper Nonlinear differential-difference equations and Fourier
analysis
\jour J. Math. Phys.\vol 17\yr 1976\pages 1011-1018
\endref

\ref\key{CMV}
\by Cantero, M., Moral, L. and Vel\'azquez, L.
\paper Five-diagonal matrices and zeros of orthogonal polynomials
on the unit circle
\jour Lin. Alg. Appl.\vol362\yr 2003\pages 29-56
\endref

\ref\key{D}
\by Drinfel'd, V.
\paper Hamiltonian structures on Lie groups, Lie bialgebra, and the
geometric meaning of the classical Yang-Baxter equations
\jour Soviet Math. Dokl.\vol 27\yr 1983 \pages 68-71
\endref

\ref\key{DLT}
\by Deift, P., Li, L.C. and Tomei, T.
\paper Toda flows with infinitely many variables
\jour J. Funct. Anal.\vol 64\yr 1985\pages 358-402
\endref

\ref\key{K}
\by Kostant, B.
\paper Quantization and representation theory
\inbook Representation theory of Lie groups, Proc. SRC/LMS Res. Symp.,
Oxford 1977, LMS Lecture Notes Series 34
\eds Atiyah, M.
\publ Cambridge University Press \yr 1979\pages 287-316
\endref

\ref\key{LW} 
\by Lu, J.-H., Weinstein, A.
\paper Poisson Lie groups, dressing transformations, and Bruhat
decompositions.
\jour J. Diff. Geom. 31\yr 1990 \pages 501--526
\endref

\ref\key{N1}
\by Nenciu, I.
\paper Lax pairs for the Ablowitz-Ladik system via orthogonal 
polynomials on the unit circle
\jour Int. Math. Res. Not. 2005\issue 11\pages 647-686
\endref

\ref\key{N2}
\by Nenciu, I.
\paper Thesis
\jour Caltech\yr 2005
\endref

\ref\key{R}
\by Ratiu, T.
\paper On the smoothness of the time $t$-map of the KdV equation and the 
bifurcation of the eigenvalues of Hill's operator
\inbook Global Analysis
\bookinfo Lecture Notes in Math.\vol 755
\publ Springer
\publaddr Berlin
\yr 1979\pages 248-294
\endref

\ref\key{S}
\by Simon, B.
\book Orthogonal polynomials on the unit circle, Parts 1 and 2
\bookinfo American Mathematical Society Colloquium Publications
\publ American Mathematical Society
\publaddr Providence, R. I.
\yr 2005
\endref

\ref\key{STS1}
\by Semenov-Tian-Shansky, M.
\paper What is a classical r-matrix?
\jour Funct. Anal. Appl.\vol 17\yr 1983\pages 259-272
\endref

\ref\key{STS2}
\by Semenov-Tian-Shansky, M.
\paper Dressing transformations and Poisson group actions
\jour Publ. RIMS, Kyoto University \vol 21\yr 1985 \pages1237-1260
\endref

\ref\key{T}
\by Teschl, G.
\book Jacobi operators and completely integrable nonlinear lattices
\bookinfo Mathematical Surveys and Monographs\vol 72
\publ American Mathematical Society
\publaddr Providence, R. I.
\yr 2000
\endref

\endRefs
\enddocument